# ISS WITH RESPECT TO BOUNDARY DISTURBANCES FOR 1-D PARABOLIC PDES

**Iasson Karafyllis[*] and Miroslav Krstic[**]**

[*]Dept. of Mathematics, National Technical University of Athens, Zografou Campus, 15780, Athens, Greece, email: iasonkar@central.ntua.gr

[**]Dept. of Mechanical and Aerospace Eng., University of California, San Diego, La Jolla, CA 92093-0411, U.S.A., email: krstic@ucsd.edu

**Abstract**

Due to unbounded input operators in partial differential equations (PDEs) with boundary inputs, there has been a long-held intuition that input-to-state stability (ISS) properties and finite gains cannot be established with respect to disturbances at the boundary. This intuition has been reinforced by many unsuccessful attempts, as well as by the success in establishing ISS only with respect to the derivative of the disturbance. Contrary to this intuition, we establish such a result for parabolic PDEs. Our methodology does not rely on the transformation of the boundary disturbance to a distributed input and the stability analysis is performed in time-varying subsets of the state space. The obtained results are used for the comparison of the gain coefficients of transport PDEs with respect to inlet disturbances and for the establishment of the ISS property with respect to control actuator errors for parabolic systems under boundary feedback control.

## 1. Introduction

The extension of the Input-to-State Stability (ISS) property to systems which are described by Partial Differential Equations (PDEs) is a challenge. Many works have recently studied possible extensions of the ISS property to PDE systems (see for example [1,2,4,5,6,9,11,12,14]). In particular, for PDE systems there are two "places" where a disturbance can appear: the domain (a distributed disturbance appearing in the partial differential equations) and the boundary (a disturbance that appears in the boundary conditions). Most of the existing results in the literature are dealing with distributed disturbances in the domain (an exception is the work [1]).

Boundary disturbances can be cast as distributed disturbances acting on the domain by means of standard transformation arguments. However, when a boundary disturbance is expressed by means of a distributed disturbance then the effect of the boundary disturbance is described by means of an unbounded operator (see the relevant discussion in [11] for inputs in infinite-dimensional systems that are expressed by means of unbounded linear operators). The advantage of the methodology is that the "disturbed problem" becomes a standard evolution equation (with inputs) in an appropriate complete linear space $X$, so that all existing tools for evolution equations can be used (e.g., semigroup of linear operators). However, such a methodology will always end up not showing the ISS property with respect to the boundary disturbance but the ISS property with respect to the boundary disturbance and some of its time derivatives (see for example [1]).

The present work is devoted to the presentation of a different methodology for studying ISS with respect to boundary disturbances. The transformation of the boundary disturbance to a domain disturbance is avoided and the effect of the disturbance is not expressed by means of an unbounded operator. The effect of the boundary disturbance $d(t)$ is expressed by a change in the state space itself: the state space is different for every time instant $t \geq 0$ and depends on the value of the disturbance. Therefore, the evolution of the state takes place in a parameterized convex set $X_{d(t)}$. The focus of the present work is on 1-D parabolic PDEs, although the methodology can be extended to other classes of PDEs as well. The proof relies on the establishment of estimates for the magnitude of certain generalized Fourier coefficients. The estimate of the appropriate weighted $L^2$ norm is obtained by the Parseval's identity.

Another difference between the present work and existing works on the ISS of PDE systems is that most of the existing works on the ISS of PDE systems are using Lyapunov functionals (see for instance [10,12,14]), while the present work does not use a Lyapunov functional. The difference is important and is a consequence of the fact that in this work the evolution of the state takes place in a parameterized convex set $X_{d(t)}$. The reader should not misunderstand the statement: we are not claiming that it is impossible to find an ISS-Lyapunov functional for a boundary disturbance. However, it is difficult to find an ISS-Lyapunov functional for boundary disturbances because the state space is different at each time instant.

The results of the present work have direct consequences to various research directions. The comparison of the gains of boundary disturbances for transport PDEs is performed and the effects of the diffusion and the boundary condition at the exit of the transport device are illustrated. The ISS property of the closed-loop system with respect to control actuator errors for backstepping boundary feedback design methodologies (see [17]) is also studied in this work.

The structure of the present work is as follows: Section 2 is devoted to the presentation of the problem and the statement of the main result (Theorem 2.2). The proof of the main result is provided in Section 3, where additional results are stated and utilized. The applications of the obtained results to transport PDEs and to the study of parabolic systems under boundary feedback control are shown in Section 4. The concluding remarks are provided in Section 5. Finally, the Appendix contains the proofs of all auxiliary results used in Section 3.

**Notation.** Throughout this paper, we adopt the following notations:
* $\Re_+ := [0, +\infty)$.
* Let $U \subseteq \Re^n$ be a set with non-empty interior and let $\Omega \subseteq \Re$ be a set. By $C^0(U; \Omega)$, we denote the class of continuous mappings on $U$, which take values in $\Omega$. By $C^k(U; \Omega)$, where $k \geq 1$, we denote the class of continuous functions on $U$, which have continuous derivatives of order $k$ on $U$ and take values in $\Omega$.
* Let $r \in C^0([0,1]; (0,+\infty))$ be given. $L_r^2([0,1])$ denotes the equivalence class of measurable functions $f : [0,1] \to \Re$ for which $\|f\|_r = \left( \int_0^1 r(z) |f(z)|^2 dz \right)^{1/2} < +\infty$. $L_r^2([0,1])$ is a Hilbert space with inner product $\langle f, g \rangle = \int_0^1 r(z) f(z) g(z) dz$.
* Let $x \in C^0(\Re_+ \times [0,1]; \Re)$ be given. We use the notation $x[t]$ to denote the profile at certain $t \geq 0$, i.e., $(x[t])(z) = x(t,z)$ for all $z \in [0,1]$.

## 2. Problem Description and Main Result

Consider the Sturm-Liouville operator $A: D \to C^0([0,1]; \Re)$ defined by

$$(Af)(z) = -\frac{1}{r(z)} \frac{d}{dz}\left(p(z) \frac{df}{dz}(z)\right) + \frac{q(z)}{r(z)} f(z), \text{ for all } f \in D \text{ and } z \in (0,1) \quad (2.1)$$

where $p \in C^1([0,1]; (0,+\infty))$, $r \in C^0([0,1]; (0,+\infty))$, $q \in C^0([0,1]; \Re)$ and $D \subseteq C^2([0,1]; \Re)$ is the set of all functions $f : [0,1] \to \Re$ for which

$$b_1 f(0) + b_2 \frac{df}{dz}(0) = a_1 f(1) + a_2 \frac{df}{dz}(1) = 0 \quad (2.2)$$

where $a_1, a_2, b_1, b_2$ are real constants with $|a_1| + |a_2| > 0$, $|b_1| + |b_2| > 0$.

**FACT** (see Chapter 11 in [3] and pages 498-505 in [13]): All eigenvalues of the Sturm-Liouville operator $A: D \to C^0([0,1]; \Re)$, defined by (2.1), (2.2) are real. They form an infinite, increasing sequence $\lambda_1 < \lambda_2 < ... < \lambda_n < ...$ with $\lim_{n \to \infty}(\lambda_n) = +\infty$. To each eigenvalue $\lambda_n \in \Re$ ($n = 1,2,...$) corresponds exactly one eigenfunction $\phi_n \in C^2([0,1]; \Re)$ that satisfies $A\phi_n = \lambda_n \phi_n$ and $b_1 \phi_n(0) + b_2 \frac{d\phi_n}{dz}(0) = a_1 \phi_n(1) + a_2 \frac{d\phi_n}{dz}(1) = 0$. The eigenfunctions form an orthonormal basis of $L_r^2([0,1])$.

In the present work, we make the following assumption for the Sturm-Liouville operator $A: D \to C^0([0,1]; \Re)$ defined by (2.1), (2.2), where $a_1, a_2, b_1, b_2$ are real constants with $|a_1| + |a_2| > 0$, $|b_1| + |b_2| > 0$.

**(H):** *The Sturm-Liouville operator $A: D \to C^0([0,1]; \Re)$ defined by (2.1), (2.2), where $a_1, a_2, b_1, b_2$ are real constants with $|a_1| + |a_2| > 0$, $|b_1| + |b_2| > 0$, satisfies*

$$\lambda_1 > 0 \quad (2.3)$$

*and*

$$\sum_{n=1}^{\infty} \lambda_n^{-1} \max_{0 \leq z \leq 1}(|\phi_n(z)|) < +\infty \quad (2.4)$$

Consider next the parameterized convex set

$$X_\mu = \left\{ x \in C^2([0,1]; \Re) : b_1 x(0) + b_2 \frac{dx}{dz}(0) = \mu, a_1 x(1) + a_2 \frac{dx}{dz}(1) = 0 \right\} \quad (2.5)$$

with parameter $\mu \in \Re$. Given $d \in C^2(\Re_+; \Re)$ and $x_0 \in X_{d(0)}$, we study the solution $x \in C^0(\Re_+ \times [0,1]; \Re) \cap C^1((0,+\infty) \times [0,1]; \Re)$ for which $x[t] \in X_{d(t)}$ for all $t \geq 0$, $x(0,z) = x_0(z)$ for all $z \in [0,1]$ and

$$\frac{\partial x}{\partial t}(t,z) + (Ax[t])(z) = \frac{\partial x}{\partial t}(t,z) - \frac{1}{r(z)} \frac{\partial}{\partial z}\left(p(z) \frac{\partial x}{\partial z}(t,z)\right) + \frac{q(z)}{r(z)} x(t,z) = 0, \text{ for all } (t,z) \in (0,+\infty) \times (0,1) \quad (2.6)$$

In other words, we consider the solution of the 1-D parabolic PDE (2.6) that satisfies for all $t \geq 0$ the boundary conditions

$$b_1 x(t,0) + b_2 \frac{\partial x}{\partial z}(t,0) = d(t) \quad , \quad a_1 x(t,1) + a_2 \frac{\partial x}{\partial z}(t,1) = 0 \tag{2.7}$$

The input $d \in C^2(\Re_+; \Re)$ is a boundary disturbance and appears only at the boundary condition.

In order to be able to state the main result, we first need the following lemma. Its proof is provided at the Appendix.

**Lemma 2.1:** *Consider the Sturm-Liouville operator $A: D \to C^0([0,1]; \Re)$ defined by (2.1), (2.2), where $a_1, a_2, b_1, b_2$ are real constants with $|a_1| + |a_2| > 0$, $b_1^2 + b_2^2 = 1$, under Assumption (H). Then the boundary value problem*

$$\frac{d}{dz}\left(p(z)\frac{d\tilde{x}}{dz}(z)\right) - q(z)\tilde{x}(z) = 0, \text{ for all } z \in [0,1], \tag{2.8}$$

*with*

$$b_1 \tilde{x}(0) + b_2 \frac{d\tilde{x}}{dz}(0) = 1 \quad , \quad a_1 \tilde{x}(1) + a_2 \frac{d\tilde{x}}{dz}(1) = 0 \tag{2.9}$$

*has a unique solution $\tilde{x} \in C^2([0,1]; \Re)$, which satisfies*

$$p^2(0) \sum_{n=1}^{\infty} \lambda_n^{-2} \left| b_1 \frac{d\phi_n}{dz}(0) - b_2 \phi_n(0) \right|^2 = \int_0^1 r(z)\tilde{x}^2(z) dz \tag{2.10}$$

We are now ready to state the main result of the present work.

**Theorem 2.2:** *Consider the Sturm-Liouville operator $A: D \to C^0([0,1]; \Re)$ defined by (2.1), (2.2), where $a_1, a_2, b_1, b_2$ are real constants with $|a_1| + |a_2| > 0$, $|b_1| + |b_2| > 0$, under assumption (H). Then for every $d \in C^2(\Re_+; \Re)$ and $x_0 \in X_{d(0)}$, the evolution equation (2.6) with (2.7) and initial condition $x_0 \in X_{d(0)}$ has a unique solution $x \in C^0(\Re_+ \times [0,1]; \Re) \cap C^1((0,+\infty) \times [0,1]; \Re)$ for which $x[t] \in X_{d(t)}$ for all $t \geq 0$, $x(0,z) = x_0(z)$ for all $z \in [0,1]$ and satisfies the following estimate for all $t \geq 0$ and $\varepsilon > 0$:*

$$\|x[t]\|_r \leq \sqrt{1+\varepsilon} \exp(-\lambda_1 t) \|x[0]\|_r + C \sqrt{\frac{1+\varepsilon^{-1}}{b_1^2 + b_2^2}} \max_{0 \leq s \leq t}(|d(s)|) \tag{2.11}$$

*where*

$$C := \frac{p(0)}{\sqrt{b_1^2 + b_2^2}} \left( \sqrt{\sum_{n=1}^{\infty} \lambda_n^{-2} \left| b_1 \frac{d\phi_n}{dz}(0) - b_2 \phi_n(0) \right|^2} \right) = \sqrt{\int_0^1 r(z)\tilde{x}^2(z) dz} \tag{2.12}$$

*and $\tilde{x} \in C^2([0,1]; \Re)$ is the unique solution of the boundary value problem (2.8) with $b_1 \tilde{x}(0) + b_2 \frac{d\tilde{x}}{dz}(0) = \sqrt{b_1^2 + b_2^2}$ and $a_1 \tilde{x}(1) + a_2 \frac{d\tilde{x}}{dz}(1) = 0$. In other words, the system described by the evolution equation (2.6) with (2.7), state space the normed linear space $\bigcup_{\mu \in \Re} X_\mu \subseteq L_r^2([0,1])$ with norm $\| \|_r$, satisfies the ISS property with respect to the boundary input $d \in C^2(\Re_+; \Re)$ with gain $C\sqrt{\frac{1+\varepsilon^{-1}}{b_1^2 + b_2^2}}$ for every $\varepsilon > 0$.*

**Remark 2.3:** Since the equilibrium point that corresponds to the constant disturbance $d(t) \equiv \sqrt{b_1^2 + b_2^2}$ is the function $\tilde{x} \in C^2([0,1]; \Re)$, it follows that the gain $\gamma > 0$ that is involved to the ISS estimate

$$\|x[t]\|_r \le M \exp(-\sigma t) \|x[0]\|_r + \gamma \max_{0 \le s \le t}(|d(s)|), \text{ for all } t \ge 0$$

for certain constants $M, \sigma$, must satisfy the inequality

$$\gamma \sqrt{b_1^2 + b_2^2} \ge C$$

where $C > 0$ is given by (2.12). On the other hand, Theorem 2.2 guarantees that

$$\gamma \le C \sqrt{\frac{1+\varepsilon^{-1}}{b_1^2 + b_2^2}} \text{ , for all } \varepsilon > 0$$

Consequently, we can guarantee that the estimation of the gain made by Theorem 2.2 is not conservative. Moreover, formula (2.12) guarantees that the gain of the boundary disturbance can be computed **without exact knowledge** of the eigenvalues and the eigenfunctions of the Sturm-Liouville operator $A: D \to C^0([0,1]; \Re)$ defined by (2.1), (2.2). The only thing we need to know about the eigenvalues and the eigenfunctions of the Sturm-Liouville operator $A: D \to C^0([0,1]; \Re)$ defined by (2.1), (2.2) is that Assumption (H) holds.

**Remark 2.4:** The ISS property guaranteed by estimate (2.11) is a direct extension of the ISS property for finite-dimensional systems (see [16]). However, there are some differences with the finite-dimensional case:
  i) The state $x[t] \in X_{d(t)}$ is not allowed to take values in the whole state space $\cup_{\mu \in \Re} X_\mu \subseteq L^2_r([0,1])$ but is restricted to evolve in the convex set $X_{d(t)}$.
  ii) The disturbance $d \in C^2(\Re_+; \Re)$ has to be sufficiently regular.

Both requirements are necessary due to the regularity requirements for the solution. Indeed, if we studied weak solutions (instead of classical solutions) of the PDE problem (2.6), (2.7), then less demanding regularity properties for the disturbance would be needed.

## 3. Proof of Main Result

In order to prove the main result, we first need an existence/uniqueness result. Although the following result guarantees the existence/uniqueness of a classical solution for a PDE problem, we have not been able to find such a result in the literature. The results in [7,14] could be applied in principle, but we would have obtained a less smooth solution. Therefore, we are forced to prove the following result. Its proof is provided at the Appendix.

**Theorem 3.1:** *Consider the Sturm-Liouville operator $A: D \to C^0([0,1]; \Re)$ defined by (2.1), (2.2), where $a_1, a_2, b_1, b_2$ are real constants with $|a_1| + |a_2| > 0$, $|b_1| + |b_2| > 0$, under Assumption (H). Consider the linear space $X_0 = \left\{ x \in C^2([0,1]; \Re) : b_1 x(0) + b_2 \frac{dx}{dz}(0) = a_1 x(1) + a_2 \frac{dx}{dz}(1) = 0 \right\}$. Then for every $x_0 \in X_0$ and $f \in C^1(\Re_+ \times [0,1]; \Re)$, there exists a unique function $x \in C^0(\Re_+ \times [0,1]; \Re) \cap C^1((0, +\infty) \times [0,1]; \Re)$ satisfying $x[t] \in X_0$ for all $t \ge 0$, $x(0, z) = x_0(z)$ for all $z \in [0,1]$ and*

$$\frac{\partial x}{\partial t}(t,z) - \frac{1}{r(z)} \frac{\partial}{\partial z}\left(p(z) \frac{\partial x}{\partial z}(t,z)\right) + \frac{q(z)}{r(z)} x(t,z) = f(t,z) \text{ , for all } (t,z) \in (0, +\infty) \times (0,1) \quad (3.1)$$

Using Theorem 3.1 we are in a position to guarantee existence/uniqueness of a classical solution for the PDE problem (2.6), (2.7).

**Corollary 3.2:** *Consider the Sturm-Liouville operator $A:D \to C^0([0,1];\Re)$ defined by (2.1), (2.2), where $a_1, a_2, b_1, b_2$ are real constants with $|a_1|+|a_2|>0$, $b_1^2+b_2^2=1$, under Assumption (H). Then for every $d \in C^2(\Re_+;\Re)$ and $x_0 \in X_{d(0)}$, there exists a unique function $x \in C^0(\Re_+ \times [0,1];\Re) \cap C^1((0,+\infty) \times [0,1];\Re)$ for which $x[t] \in X_{d(t)}$ for all $t \geq 0$, $x(0,z) = x_0(z)$ for all $z \in [0,1]$ and (2.6).*

**Proof:** We simply apply the transformation $x(t,z) = y(t,z) + d(t)\big(b_1 + b_2 z + c_1 z^2 + c_2 z^3\big)$, where $y \in C^0(\Re_+ \times [0,1];\Re) \cap C^1((0,+\infty) \times [0,1];\Re)$ is the unique function that satisfies $y[t] \in X_0$ for all $t \geq 0$, $y(0,z) = x_0(z) - d(0)\big(b_1 + b_2 z + c_1 z^2 + c_2 z^3\big)$ for all $z \in [0,1]$ and

$$\frac{\partial y}{\partial t}(t,z) - \frac{1}{r(z)} \frac{\partial}{\partial z}\left(p(z) \frac{\partial y}{\partial z}(t,z)\right) + \frac{q(z)}{r(z)} y(t,z) =$$
$$\frac{d(t)}{r(z)}\left(\frac{d}{dz}\big((b_2 + 2c_1 z + 3c_2 z^2)p(z)\big) - q(z)(b_1 + b_2 z + c_1 z^2 + c_2 z^3)\right) - \dot{d}(t)(b_1 + b_2 z + c_1 z^2 + c_2 z^3),$$
for all $(t,z) \in (0,+\infty) \times (0,1)$ (3.2)

and $c_1, c_2 \in \Re$ are constants that satisfy $(a_1 + 2a_2)c_1 + (a_1 + 3a_2)c_2 = -a_1 b_1 - (a_1 + a_2)b_2$. Notice that the existence of constants $c_1, c_2 \in \Re$ that satisfy $(a_1 + 2a_2)c_1 + (a_1 + 3a_2)c_2 = -a_1 b_1 - (a_1 + a_2)b_2$ is guaranteed by the condition $|a_1|+|a_2|>0$. Moreover, the existence/uniqueness of $y \in C^0(\Re_+ \times [0,1];\Re) \cap C^1((0,+\infty) \times [0,1];\Re)$ that satisfies $y[t] \in X_0$ for all $t \geq 0$, $y(0,z) = x_0(z) - d(0)\big(b_1 + b_2 z + c_1 z^2 + c_2 z^3\big)$ for all $z \in [0,1]$ and (3.2) is guaranteed by Theorem 3.1. The proof is complete. ◁

We are now ready to prove Theorem 2.1.

**Proof of Theorem 2.1:** We first restrict our attention to the case where $b_1^2 + b_2^2 = 1$. Corollary 3.2 guarantees that $x[t] \in X_{d(t)} \subset L_r^2([0,1])$ for all $t \geq 0$. Since the eigenfunctions $\{\phi_n\}_{n=1}^{\infty}$ of the Sturm-Liouville operator $A:D \to C^0([0,1];\Re)$ defined by (2.1), (2.2) form an orthonormal basis of $L_r^2([0,1])$, it follows that Parseval's identity holds, i.e.,

$$\|x[t]\|_r^2 = \sum_{n=1}^{\infty} c_n^2(t), \text{ for all } t \geq 0 \tag{3.3}$$

where

$$c_n(t) := \langle \phi_n, x[t] \rangle = \int_0^1 r(z) x(t,z) \phi_n(z) dz, \text{ for } n = 1, 2, \ldots \tag{3.4}$$

By virtue of (2.6), (2.7) and the facts that $-\frac{1}{r(z)} \frac{d}{dz}\left(p(z) \frac{d\phi_n}{dz}(z)\right) + \frac{q(z)}{r(z)} \phi_n(z) = \lambda_n \phi_n(z)$, $b_1 \phi_n(0) + b_2 \frac{d\phi_n}{dz}(0) = a_1 \phi_n(1) + a_2 \frac{d\phi_n}{dz}(1) = 0$, it follows from repeated integration by parts, that the following equalities hold for all $t > 0$:

$$\dot{c}_n(t) = \int_0^1 r(z) \frac{\partial x}{\partial t}(t,z) \phi_n(z) dz$$

$$= \int_0^1 \frac{\partial}{\partial z}\left(p(z) \frac{\partial x}{\partial z}(t,z)\right) \phi_n(z) dz - \int_0^1 q(z) x(t,z) \phi_n(z) dz$$

$$= p(1) \frac{\partial x}{\partial z}(t,1) \phi_n(1) - p(0) \frac{\partial x}{\partial z}(t,0) \phi_n(0) - \int_0^1 p(z) \frac{\partial x}{\partial z}(t,z) \frac{d\phi_n}{dz}(z) dz - \int_0^1 q(z) x(t,z) \phi_n(z) dz =$$

$$= p(1)\left(\frac{\partial x}{\partial z}(t,1) \phi_n(1) - x(t,1) \frac{d\phi_n}{dz}(1)\right) + p(0)\left(\frac{d\phi_n}{dz}(0) x(t,0) - \frac{\partial x}{\partial z}(t,0) \phi_n(0)\right) + \int_0^1 x(t,z)\left[\frac{d}{dz}\left(p(z) \frac{d\phi_n}{dz}(z)\right) - q(z) \phi_n(z)\right] dz$$

$$= p(1)\left(\frac{\partial x}{\partial z}(t,1) \phi_n(1) - x(t,1) \frac{d\phi_n}{dz}(1)\right) + p(0)\left(\frac{d\phi_n}{dz}(0) x(t,0) - \frac{\partial x}{\partial z}(t,0) \phi_n(0)\right) - \int_0^1 r(z) x(t,z) (A\phi_n)(z) dz$$

$$= p(0)\left(b_1 \frac{d\phi_n}{dz}(0) - b_2 \phi_n(0)\right) d(t) - \lambda_n c_n(t)$$

In the above equations, we have used the fact that by virtue of (2.7) and $a_1 \phi_n(1) + a_2 \frac{d\phi_n}{dz}(1) = 0$, the homogeneous system of linear equations

$$s_1 x(t,1) + s_2 \frac{\partial x}{\partial z}(t,1) = 0$$

$$s_1 \phi_n(1) + s_2 \frac{d\phi_n}{dz}(1) = 0$$

has the non-zero solution

$$s_1 = a_1 \quad , \quad s_2 = a_2$$

and consequently, the determinant of the matrix $\begin{bmatrix} x(t,1) & \frac{\partial x}{\partial z}(t,1) \\ \phi_n(1) & \frac{d\phi_n}{dz}(1) \end{bmatrix}$ is zero, i.e.,

$x(t,1) \frac{d\phi_n}{dz}(1) - \phi_n(1) \frac{\partial x}{\partial z}(t,1) = 0$. Moreover, we have used the fact that

$x(t,0) \frac{d\phi_n}{dz}(0) - \frac{\partial x}{\partial z}(t,0) \phi_n(0) = d(t)\left(b_1 \frac{d\phi_n}{dz}(0) - b_2 \phi_n(0)\right)$, which is a direct consequence of the facts

$b_1 x(t,0) + b_2 \frac{\partial x}{\partial z}(t,0) = d(t)$, $b_1 \phi_n(0) + b_2 \frac{d\phi_n}{dz}(0) = 0$ and $b_1^2 + b_2^2 = 1$.

Integrating the above differential equations, we obtain for all $0 < T \leq t$ and $n = 1, 2, \ldots$:

$$c_n(t) = \exp(-\lambda_n (t-T)) c_n(T) - p(0)\left(b_1 \frac{d\phi_n}{dz}(0) - b_2 \phi_n(0)\right) \int_T^t \exp(-\lambda_n (t-s)) d(s) ds \tag{3.5}$$

Continuity of the mapping $\Re_+ \ni T \to c_n(T)$ and (3.5) implies the following equations for all $t \geq 0$ and $n = 1, 2, \ldots$:

$$c_n(t) = \exp(-\lambda_n t) c_n(0) - p(0)\left(b_1 \frac{d\phi_n}{dz}(0) - b_2 \phi_n(0)\right) \int_0^t \exp(-\lambda_n (t-s)) d(s) ds \tag{3.6}$$

Equations (3.6) imply the following estimates for all $t \geq 0$ and $n = 1, 2, \ldots$:

$$\left|c_n(t)\right| \leq \exp(-\lambda_n t)\left|c_n(0)\right| + p(0)\left|b_1 \frac{d\phi_n}{dz}(0) - b_2\phi_n(0)\right| \frac{1 - \exp(-\lambda_n t)}{\lambda_n} \max_{0 \leq s \leq t}\left(\left|d(s)\right|\right) \quad (3.7)$$

Since $\lambda_n \geq \lambda_1 > 0$ for all $n = 1, 2, \ldots$, we obtain the following estimates for all $t \geq 0$, $\varepsilon > 0$ and $n = 1, 2, \ldots$:

$$\left|c_n(t)\right|^2 \leq (1+\varepsilon)\exp(-2\lambda_1 t)\left|c_n(0)\right|^2 + (1+\varepsilon^{-1})p^2(0)\left|b_1 \frac{d\phi_n}{dz}(0) - b_2\phi_n(0)\right|^2 \frac{1}{\lambda_n^2} \max_{0 \leq s \leq t}\left(\left|d(s)\right|^2\right) \quad (3.8)$$

Therefore, by virtue of estimates (3.8), (3.3), definition (2.12) and identity (2.10), the following estimate holds for all $t \geq 0$ and $\varepsilon > 0$:

$$\|x[t]\|_r \leq \sqrt{1+\varepsilon} \exp(-\lambda_1 t)\|x[0]\|_r + C\sqrt{1+\varepsilon^{-1}} \max_{0 \leq s \leq t}\left(\left|d(s)\right|\right) \quad (3.9)$$

where $C > 0$ is given by (2.12).

Next, consider the general case, where the constants $b_1, b_2$ appearing in the boundary condition $b_1 x(t,0) + b_2 \frac{\partial x}{\partial z}(t,0) = d(t)$ does not satisfy the condition $b_1^2 + b_2^2 = 1$. Notice that the boundary condition $b_1 x(t,0) + b_2 \frac{\partial x}{\partial z}(t,0) = d(t)$ can be transformed to the condition $\tilde{b}_1 x(t,0) + \tilde{b}_2 \frac{\partial x}{\partial z}(t,0) = \tilde{d}(t)$, where $\tilde{b}_i = \frac{b_i}{\sqrt{b_1^2 + b_2^2}}$ ($i = 1,2$) and $\tilde{d}(t) = \frac{d(t)}{\sqrt{b_1^2 + b_2^2}}$. Therefore, (3.9) holds with $d$ replaced by $\tilde{d}$. Estimate (2.11) is a direct consequence of estimate (3.9) with $d$ replaced by $\tilde{d}$ and definition $\tilde{d}(t) = \frac{d(t)}{\sqrt{b_1^2 + b_2^2}}$. The proof is complete. ◁

**Remark 3.3:** It is clear that the proof of Theorem 2.1 consists of two different parts. The first part that establishes existence/uniqueness relies on the transformation of the boundary disturbance to a domain input. However, this transformation is used only for the establishment of existence/uniqueness of the solution. Estimates of the magnitude of certain generalized Fourier coefficients are derived in the second part: the estimates are derived for the original equations not the transformed ones. The obtained estimates are combined by means of Parseval's identity in order to give an estimate for the $L^2$ norm.

## 4. Applications

**4.I. Gains for Transport PDEs With Respect to Inlet Disturbances**

We consider the 1-D transport PDE

$$\frac{\partial y}{\partial t}(t,z) = D\frac{\partial^2 y}{\partial z^2}(t,z) - v\frac{\partial y}{\partial z}(t,z) - ky(t,z) \quad (4.1)$$

where $D > 0$, $v \geq 0$, $k \in \Re$ are constants. We consider the following cases.

CASE 1: Dirichlet boundary conditions
$$\begin{aligned} y(t,0) &= d(t) \\ y(t,1) &= 0 \end{aligned} \quad (4.2)$$

CASE 2: Robin (or Neumann) boundary conditions

$$y(t,0) = d(t)$$
$$\frac{\partial y}{\partial z}(t,1) = \left(\frac{v}{2D} - a\right) y(t,1) \quad (4.3)$$

where $a \geq 0$.

The physical meaning of the 1-D transport PDE is the transportation of a certain quantity (denoted by $y$) through a tube (the transport device). The inlet of the tube is at $z = 0$ while the exit is at $z = 1$. The term $D\frac{\partial^2 y}{\partial z^2}(t,z)$ quantifies the effect of the diffusion and $D > 0$ is called the diffusion coefficient, the term $-v\frac{\partial y}{\partial z}(t,z)$ quantifies the effect of convection and $v \geq 0$ is the fluid velocity in the tube, while the term $-ky(t,z)$ quantifies the possible reaction effects and $k \in \Re$ is the reaction constant. Notice that we are considering a parameterized family of boundary conditions at the exit of the tube with parameter $a \geq 0$ (Case 1 corresponds to $a = +\infty$). The disturbance is at the inlet of the tube and is transported throughout the tube by means of diffusion and convection.

The 1-D transport PDE (4.1) corresponds to the PDE (2.6) with

$$p(z) = D\exp\left(-\frac{v}{D}z\right) \ , \ r(z) = \exp\left(-\frac{v}{D}z\right) \ , \ q(z) = k\exp\left(-\frac{v}{D}z\right) \quad (4.4)$$

At this point, we could proceed to the analysis of the 1-D transport PDE (4.1) with boundary condition given either by (4.2) (which corresponds to the case $a = +\infty$) or by (4.3). However, in order to make the algebraic manipulations easier, we apply the transformation $y(t,z) = \exp\left(\frac{vz}{2D}\right) x(t,z)$, which transforms the PDE (4.1) to the following PDE

$$\frac{\partial x}{\partial t}(t,z) = D\frac{\partial^2 x}{\partial z^2}(t,z) - \left(\frac{v^2}{4D} + k\right) x(t,z) \quad (4.5)$$

with the following boundary conditions

CASE 1: Dirichlet boundary conditions

$$x(t,0) = d(t)$$
$$x(t,1) = 0 \quad (4.6)$$

CASE 2: Robin (or Neumann) boundary conditions

$$x(t,0) = d(t)$$
$$\frac{\partial x}{\partial z}(t,1) = -a\, x(t,1) \quad (4.7)$$

The 1-D transport PDE (4.5) corresponds to the PDE (2.6) with

$$p(z) \equiv D \ , \ r(z) \equiv 1 \ , \ q(z) \equiv k + \frac{v^2}{4D} \quad (4.8)$$

In every case, the eigenvalues are

$$\lambda_n = k + \frac{v^2}{4D} + D\omega_n^2, \quad n = 1, 2, \ldots \tag{4.9}$$

and the eigenfunctions are

$$\phi_n(z) = \sqrt{\frac{2}{1 - \frac{\sin(2\omega_n)}{2\omega_n}}} \sin(\omega_n z), \quad n = 1, 2, \ldots \tag{4.10}$$

where $\omega_n = (n - \mu_n(a))\pi$ and $\mu_n(a) \in \left[0, \frac{1}{2}\right]$ ($n = 1, 2, \ldots$) are given by:

- $\mu_n(+\infty) = 0$, for Case 1 that corresponds to $a = +\infty$,

- $\mu_n(a) \in \left(0, \frac{1}{2}\right)$ for Case 2 with $a > 0$ is the unique solution of the equation $\tan(\mu_n \pi) + a^{-1} \mu_n \pi = a^{-1} n \pi$,

- $\mu_n(0) = \frac{1}{2}$, for Case 2 with $a = 0$.

The assumption $\lambda_1 > 0$ is equivalent to the following condition:

$$k > -\frac{v^2}{4D} - D\pi^2 (1 - \mu_1(a))^2 \tag{4.11}$$

Since $2\omega_n \geq \pi$ for every case, it follows from (4.10) that $\max_{0 \leq z \leq 1}(|\phi_n(z)|) \leq \sqrt{\frac{2\pi}{\pi - 1}}$. Moreover, since $\omega_n \geq \left(n - \frac{1}{2}\right)\pi$ for every case, it follows that $\lambda_n \geq k + \frac{v^2}{4D} + D\pi^2 \left(n - \frac{1}{2}\right)^2$. Therefore, for $N > \frac{1}{2} + \frac{1}{2D\pi}\sqrt{\max(0, -v^2 - 4kD)}$, we get:

$$\sum_{n=1}^{\infty} \lambda_n^{-1} \max_{0 \leq z \leq 1}(|\phi_n(z)|) \leq \sqrt{\frac{2\pi}{\pi - 1}} N \lambda_1^{-1} + \sqrt{\frac{2\pi}{\pi - 1}} \sum_{n=N+1}^{\infty} \frac{4D}{4kD + v^2 + 4D^2\pi^2 \left(n - \frac{1}{2}\right)^2} < +\infty$$

Therefore, condition (2.4) holds for every case. In what follows, we will assume that

$$k > -\frac{v^2}{4D} \tag{4.12}$$

When (4.12) holds, the solution $\tilde{x} \in C^2([0,1]; \Re)$ of the boundary value problem $D\frac{d^2 \tilde{x}}{dz^2}(z) - \left(\frac{v^2}{4D} + k\right)\tilde{x}(z) = 0$ with $\tilde{x}(0) = 1$ and $\frac{d\tilde{x}}{dz}(1) = -a\tilde{x}(1)$ is given by

$$\tilde{x}(z) = c_1 \exp(\zeta z) + c_2 \exp(-\zeta z), \text{ for } z \in [0, 1] \tag{4.13}$$

where $\zeta := \frac{1}{2D}\sqrt{v^2 + 4kD}$, $c_1 := \frac{\zeta - a}{(\zeta + a)\exp(2\zeta) + \zeta - a}$, $c_2 := \frac{(\zeta + a)\exp(2\zeta)}{(\zeta + a)\exp(2\zeta) + \zeta - a}$. It follows from Lemma 2.1 that

$$G(\zeta,a) := \frac{\sqrt{2}}{\pi}\sqrt{\sum_{n=1}^{\infty}\frac{(n-\mu_n(a))^3}{\left((n-\mu_n(a))+(2\pi)^{-1}\sin(2\mu_n(a)\pi)\right)\left(\pi^{-2}\zeta^2+(n-\mu_n(a))^2\right)^2}} \quad (4.14)$$

$$= \sqrt{c_1^2\frac{\exp(2\zeta)-1}{2\zeta}+c_2^2\frac{1-\exp(-2\zeta)}{2\zeta}+2c_1 c_2}$$

When (4.12) holds, the solution $\tilde{x} \in C^2([0,1];\Re)$ of the boundary value problem $D\frac{d^2\tilde{x}}{dz^2}(z)-\left(\frac{v^2}{4D}+k\right)\tilde{x}(z)=0$ with $\tilde{x}(0)=1$ and $\tilde{x}(1)=0$ is given by (4.13), where $\zeta := \frac{1}{2D}\sqrt{v^2+4kD}$, $c_1 := \frac{-1}{\exp(2\zeta)-1}$, $c_2 := \frac{\exp(2\zeta)}{\exp(2\zeta)-1}$. It follows from Lemma 2.1 that

$$G(\zeta,+\infty) := \frac{\sqrt{2}}{\pi}\sqrt{\sum_{n=1}^{\infty}\frac{n^2}{\left(\pi^{-2}\zeta^2+n^2\right)^2}} = \frac{1}{\exp(2\zeta)-1}\sqrt{\frac{\exp(4\zeta)-1-4\zeta\exp(2\zeta)}{2\zeta}} \quad (4.15)$$

We are now ready to apply Theorem 2.2 and take into account the transformation $y(t,z) = \exp\left(\frac{vz}{2D}\right)x(t,z)$.

**Corollary 4.1:** *Consider the PDE (4.1) with boundary given either by (4.2) (which corresponds to $a=+\infty$) or by (4.3) with $a \geq 0$. Suppose that (4.12) holds. Then for every $d \in C^2(\Re_+;\Re)$ and $y_0 \in C^2([0,1];\Re)$ with $y_0(0)=d(0)$ and $y_0(1)=0$ (in case that (4.2) holds) or $\frac{dy_0}{dz}(1)=\left(\frac{v}{2D}-a\right)y_0(1)$ (in case that (4.3) holds), the solution of the evolution equation (4.1) with either (4.2) or (4.3) and initial condition $y_0$ is unique, defined for all $t \geq 0$ and satisfies the following estimate for all $t \geq 0$ and $\varepsilon > 0$:*

$$\sqrt{\int_0^1 \exp\left(-\frac{vz}{D}\right)y^2(t,z)dz}$$
$$\leq \exp\left(-D\left(\zeta^2+\pi^2(1-\mu_1(a))^2\right)t\right)\sqrt{(1+\varepsilon)\int_0^1 \exp\left(-\frac{vz}{D}\right)y_0^2(z)dz} + \sqrt{1+\varepsilon^{-1}}\,G(\zeta,a)\max_{0\leq s\leq t}(|d(s)|) \quad (4.16)$$

where $\zeta := \frac{1}{2D}\sqrt{v^2+4kD}$.

Next consider, for comparison purposes, the advection equation

$$\frac{\partial y}{\partial t}(t,z)+v\frac{\partial y}{\partial z}(t,z) = -ky(t,z),\text{ for all } t>0,\ z\in(0,1) \quad (4.17)$$

with $v>0$, $k>-\frac{v^2}{4D}$ and boundary condition

$$y(t,0)=d(t)\ ,\text{ for all } t>0 \quad (4.18)$$

where $d \in C^1(\Re_+;\Re)$ is a given function. The PDE (4.17) is accompanied by the initial condition

$$y(0, z) = y_0(z) \text{ for all } z \in [0,1] \tag{4.19}$$

where $y_0 \in C^1([0,1]; \Re)$ is the initial condition that satisfies $y_0(0) = d(0)$ and $\dot{d}(0) = -v \frac{dy_0}{dz}(0) - kd(0)$. The solution of (4.17), (4.18), (4.19) is given by the formula:

$$y(t,z) = (y[t])(z) = \begin{cases} \exp(-kt) y_0(z - vt) & \text{for } vt < z \\ \exp(-kv^{-1}z) d(t - v^{-1}z) & \text{for } vt \geq z \end{cases}, \text{ for all } t \geq 0, z \in [0,1] \tag{4.20}$$

Using (4.20) and the fact that $\zeta := \frac{1}{2D}\sqrt{v^2 + 4kD}$, we obtain the following estimate for all $t \geq 0$:

$$\sqrt{\int_0^1 \exp\left(-\frac{vz}{D}\right) y^2(t,z) dz} \leq \exp(-D(\zeta^2 + \pi^2 l^{-2})t) \sqrt{\int_0^1 \exp\left(-\frac{vz}{D}\right) y_0^2(z) dz} + \sqrt{\frac{1 - \exp(-(l\pi^{-1}\zeta^2 + \pi d^{-1}))}{l\pi^{-1}\zeta^2 + \pi d^{-1}}} \max_{\max(0, t-v^{-1}) \leq s \leq t} (|d(s)|) \tag{4.21}$$

where $l := \frac{2D}{v\pi^2}$. Indeed, using (4.20) we obtain when $vt < 1$:

$$\int_0^1 \exp(-vD^{-1}z) y^2(t,z) dz = \int_0^{vt} \exp(-vD^{-1}z) y^2(t,z) dz + \int_{vt}^1 \exp(-vD^{-1}z) y^2(t,z) dz$$

$$= \exp(-2kt) \int_{vt}^1 \exp(-vD^{-1}z) y_0^2(z - vt) dz + \int_0^{vt} \exp(-(vD^{-1} + 2kv^{-1})z) d^2(t - v^{-1}z) dz$$

$$\leq \exp(-(2k + v^2 D^{-1})t) \int_0^{1-vt} \exp(-vD^{-1}s) y_0^2(s) ds + \int_0^1 \exp(-(vD^{-1} + 2kv^{-1})z) dz \max_{\max(0,t-v^{-1}) \leq s \leq t} (|d(s)|^2)$$

$$\leq \exp(-(2k + v^2 D^{-1})t) \int_0^1 \exp(-vD^{-1}s) y_0^2(s) ds + \frac{1 - \exp(-(vD^{-1} + 2kv^{-1}))}{vD^{-1} + 2kv^{-1}} \max_{\max(0,t-v^{-1}) \leq s \leq t} (|d(s)|^2)$$

The same estimate is obtained when $vt \geq 1$ as well. It follows from the above estimate and definitions $\zeta := \frac{1}{2D}\sqrt{v^2 + 4kD}$, $l := \frac{2D}{v\pi^2}$ that estimate (4.21) holds.

Figure 1 shows the gains $G(\zeta, a)$ of the PDE (4.1) with respect to $\zeta$ for three different values of $a$ ($a = 0, 1, +\infty$) and the gain $\sqrt{\frac{1 - \exp(-(l\pi^{-1}\zeta^2 + \pi d^{-1}))}{l\pi^{-1}\zeta^2 + \pi d^{-1}}}$ of the advection equation (4.17) (as predicted by (4.21)). The gains have been computed under the condition $k = 0$ (no reaction), which implies that $l = \pi^{-2}\zeta^{-1}$. It is clearly shown that the gains of the PDE (4.1) are decreasing with $a$ and the case that guarantees the smallest gains is the case of Dirichlet boundary conditions ($a = +\infty$). On the other hand, the gain of the advection equation (4.17) is smaller than the gains of the PDE (4.1), except for the case that $\zeta = \frac{v}{2D}$ is small. We conclude that for the no reaction case ($k = 0$), the PDE (4.1) with diffusion has lower gain than the advection equation (4.17) only when $v \ll D$ and $a$ is large.

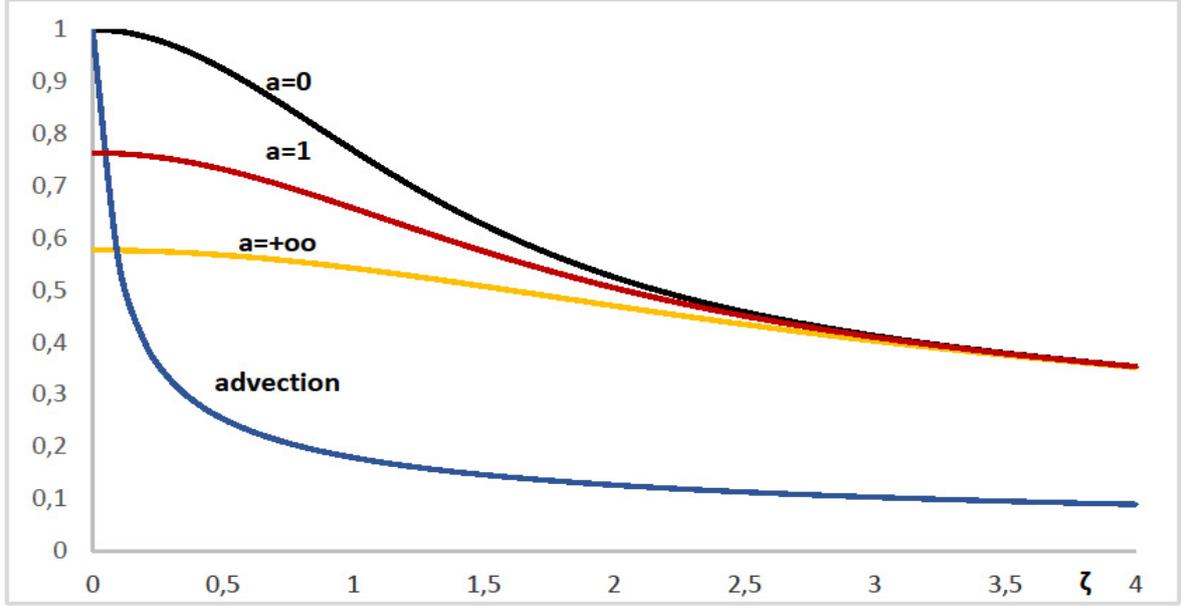

**Figure 1:** Graphs of $G(\zeta, a)$ with respect to $\zeta$ for three different values of $a$ ($a = 0, 1, +\infty$) and the graph of the gain of the advection equation (4.17). For all cases $k = 0$.

**Remark 4.2:** At this point, the reader may be surprised that the gain of the advection equation (4.17) depends on the diffusion coefficient $D > 0$ (through $\zeta := \frac{1}{2D}\sqrt{v^2 + 4kD}$ and $l := \frac{2D}{v\pi^2}$). It is clear that the solution (4.20) of the advection equation (4.17) does not depend on $D > 0$. However, the weight function $\exp\left(-\frac{vz}{D}\right)$, used in the weighted $L^2$ norm $\sqrt{\int_0^1 \exp\left(-\frac{vz}{D}\right) y^2(t,z) dz}$ that appears in both the left hand sides of (4.16) and (4.21), depends on the diffusion coefficient $D > 0$ and consequently, the gain of the advection equation (4.17) depends on the diffusion coefficient $D > 0$.

**4.II. ISS With Respect to Control Actuator Errors for Boundary State Feedback**

The recent work [17] proposed the exponential stabilization of parabolic PDEs of the form

$$\frac{\partial y}{\partial t}(t,z) = D\frac{\partial^2 y}{\partial z^2}(t,z) + py(t,z), \text{ for all } (t,z) \in (0,+\infty) \times (0,1) \quad (4.22)$$

where $D > 0$, $p \in \Re$ are constants, subject to the boundary conditions

$$\begin{array}{l} y(t,0) = u(t) \\ y(t,1) = 0 \end{array}, \text{ for all } t \geq 0 \quad (4.23)$$

where $u(t) \in \Re$ is the control input, by means of a boundary feedback stabilizer of the form

$$u(t) = -\int_0^1 k(0,s) y(t,s) ds, \text{ for all } t \geq 0 \quad (4.24)$$

where $k \in C^2([0,1]^2;\Re)$ is an appropriate function. The function $k \in C^2([0,1]^2;\Re)$ is obtained as the Volterra kernel of a Volterra integral transformation of the form

$$x(t,z) = y(t,z) + \int_z^1 k(z,s)y(t,s)ds, \text{ for all } (t,z) \in \Re_+ \times [0,1] \tag{4.25}$$

which transforms the PDE problem (4.22), (4.23), (4.24) to the problem

$$\frac{\partial x}{\partial t}(t,z) = D\frac{\partial^2 x}{\partial z^2}(t,z) - cx(t,z), \text{ for all } (t,z) \in (0,+\infty) \times (0,1) \tag{4.26}$$

where $c \geq 0$, subject to the boundary conditions

$$x(t,0) = x(t,1) = 0, \text{ for all } t \geq 0 \tag{4.27}$$

The free parameter $c \geq 0$ can be used to set the convergence rate. The solution of the original problem can be found by the inverse Volterra integral transformation

$$y(t,z) = x(t,z) + \int_z^1 l(z,s)x(t,s)ds, \text{ for all } (t,z) \in \Re_+ \times [0,1] \tag{4.28}$$

where $l \in C^2([0,1]^2;\Re)$ is an appropriate kernel.

It should be remarked that in [17] the control input is applied at $z=1$ instead of $z=0$, but the transformation of the spatial variable $z \to 1-z$ allows the statement of the results in the above form (with the control action applied at $z=0$). Moreover, it should be remarked that in [17], more general cases than the case (4.22), (4.23) are studied. Due to the similarity of all cases to the case (4.22), (4.23), we restrict our attention to the case (4.22), (4.23).

When control actuator errors are present, i.e., when the applied control action is of the form

$$u(t) = d(t) - \int_0^1 k(0,s)y(t,s)ds, \text{ for all } t \geq 0 \tag{4.29}$$

where $d \in C^2(\Re_+;\Re)$, then the transformed solution $x(t,z)$ satisfies (4.26) subject to the boundary conditions

$$\begin{aligned} x(t,0) &= d(t) \\ x(t,1) &= 0 \end{aligned}, \text{ for all } t \geq 0 \tag{4.30}$$

The PDE (4.26) corresponds to the PDE (2.6) with

$$p(z) \equiv D \quad, \quad r(z) \equiv 1 \quad, \quad q(z) \equiv c \tag{4.31}$$

The eigenvalues are

$$\lambda_n = c + Dn^2\pi^2, \quad n = 1,2,... \tag{4.32}$$

and the eigenfunctions are

$$\phi_n(z) = \sqrt{2}\sin(n\pi z), \quad n = 1,2,... \tag{4.33}$$

and assumption (H) holds. Therefore, we can apply Theorem 2.2 (exactly as in the previous section) and obtain the following estimate for the solution $x(t,z)$ of (4.26), (4.30)

$$\|x[t]\| \leq \exp\left(-\left(c+D\pi^2\right)t\right)\|x_0\|\sqrt{(1+\varepsilon)} + \sqrt{1+\varepsilon^{-1}}\, G \max_{0 \leq s \leq t}(|d(s)|), \text{ for all } t \geq 0 \text{ and } \varepsilon > 0 \quad (4.34)$$

where $\|x\| = \sqrt{\int_0^1 x^2(z)dz}$ is the standard $L^2$ norm,

$$G := \frac{\sqrt{2}}{\pi}\sqrt{\sum_{n=1}^{\infty}\frac{n^2}{\left(\pi^{-2}\zeta^2 + n^2\right)^2}} = \frac{1}{\exp(2\zeta)-1}\sqrt{\frac{\exp(4\zeta) - 1 - 4\zeta\exp(2\zeta)}{2\zeta}}, \text{ for } \zeta := \sqrt{c/D} > 0 \quad (4.35)$$

$$G := \frac{\sqrt{2}}{\pi}\sqrt{\sum_{n=1}^{\infty}\frac{1}{n^2}} = \frac{1}{\sqrt{3}}, \text{ for } c = 0 \quad (4.36)$$

for every $d \in C^2(\Re_+; \Re)$ and for every initial condition $x_0 \in C^2([0,1]; \Re)$ with $x_0(0) = d(0)$ and $x_0(1) = 0$. Using the inverse transformation (4.28), we obtain for all $t \geq 0, z \in [0,1]$ and $\varphi > 0$:

$$y^2(t,z) \leq (1+\varphi)x^2(t,z) + (1+\varphi^{-1})\left(\int_z^1 l(z,s)x(t,s)ds\right)^2 \quad (4.37)$$

Using the Cauchy-Schwarz inequality, we obtain for all $t \geq 0, z \in [0,1]$:

$$\left(\int_z^1 l(z,s)x(t,s)ds\right)^2 \leq \left(\int_z^1 l^2(z,s)ds\right)\left(\int_z^1 x^2(t,s)ds\right) \leq \left(\int_z^1 l^2(z,s)ds\right)\|x[t]\|^2 \quad (4.38)$$

Combining (4.37), (4.38) and taking $\varphi = \sqrt{\int_0^1\left(\int_z^1 l^2(z,s)ds\right)dz}$, we obtain for all $t \geq 0$:

$$\|y[t]\|^2 \leq \left(1 + \sqrt{\int_0^1\left(\int_z^1 l^2(z,s)ds\right)dz}\right)^2 \|x[t]\|^2 \quad (4.39)$$

Similarly, using (4.25) we obtain for all $t \geq 0$:

$$\|x[t]\|^2 \leq \left(1 + \sqrt{\int_0^1\left(\int_z^1 k^2(z,s)ds\right)dz}\right)^2 \|y[t]\|^2 \quad (4.40)$$

Consequently, for every $d \in C^2(\Re_+; \Re)$ and $y_0 \in C^2([0,1]; \Re)$ with $y_0(0) = d(0) - \int_0^1 k(0,s)y_0(s)ds$ and $y_0(1) = 0$, the solution of the closed-loop system (4.22), (4.23) with (4.29) and initial condition $y_0$ is unique, defined for all $t \geq 0$ and satisfies the following estimate for all $t \geq 0$ and $\varepsilon > 0$:

$$\|y[t]\| \leq \exp\left(-\left(c+D\pi^2\right)t\right)\left(1+\sqrt{\int_0^1\left(\int_z^1 l^2(z,s)ds\right)dz}\right)\left(1+\sqrt{\int_0^1\left(\int_z^1 k^2(z,s)ds\right)dz}\right)\|y_0\|\sqrt{(1+\varepsilon)}$$
$$+\left(1+\sqrt{\int_0^1\left(\int_z^1 l^2(z,s)ds\right)dz}\right)\sqrt{1+\varepsilon^{-1}}\,G\max_{0\leq s\leq t}\left(|d(s)|\right)$$
(4.41)

where $G > 0$ is defined by (4.35) and (4.36).

Estimate (4.41) shows that the closed-loop system (4.22), (4.23) with (4.29) satisfies the ISS property with respect to control actuator errors. Moreover, the estimation of the gain with respect to control actuator errors provided by estimate (4.41) is $\left(1+\sqrt{\int_0^1\left(\int_z^1 l^2(z,s)ds\right)dz}\right)\sqrt{1+\varepsilon^{-1}}\,G$, for all $\varepsilon > 0$.

## 5. Concluding Remarks

A methodology for the establishment of the Input-to-State Stability (ISS) property for 1-D parabolic Partial Differential Equations (PDEs) with boundary disturbances was proposed. The methodology does not rely on the transformation of the boundary disturbance to a domain input and the stability analysis is performed in time-varying subsets of the state space. The obtained results were used for the comparison of the gain coefficients of transport PDEs with respect to inlet disturbances and for the establishment of the ISS property with respect to control actuator errors for parabolic systems under boundary feedback control.

Future work may involve the establishment of the ISS property with the $L^\infty$ norm (instead of the $L^2$ norm that was used in the present work). Novel mathematical results will be needed for this purpose, because the analogue of Parseval's identity for the $L^\infty$ norm is not available.

# Appendix

**Proof of Lemma 2.1:** Since the Sturm-Liouville operator $A: D \to C^0([0,1]; \Re)$ defined by (2.1), (2.2) satisfies $\lambda_1 > 0$, it follows that for every $f \in C^0([0,1]; \Re)$, the solution $y \in C^2([0,1]; \Re)$ of the boundary value problem $Ay = f$ with $b_1 y(0) + b_2 \frac{dy}{dz}(0) = a_1 y(1) + a_2 \frac{dy}{dz}(1) = 0$ exists and is unique.

Let $c_1, c_2 \in \Re$ be constants that satisfy $a_1(b_1 + b_2 + c_1 + c_2) + a_2(b_2 + 2c_1 + 3c_2) = 0$. Let $\tilde{x} \in C^2([0,1]; \Re)$ be the function $\tilde{x}(z) = y(z) + g(z)$ for $z \in [0,1]$, where $g(z) = b_1 + b_2 z + c_1 z^2 + c_2 z^3$, $y \in C^2([0,1]; \Re)$ is the unique solution of the boundary value problem $Ay = -Ag$ with $b_1 y(0) + b_2 \frac{dy}{dz}(0) = a_1 y(1) + a_2 \frac{dy}{dz}(1) = 0$. It follows from the equations $b_1^2 + b_2^2 = 1$, $a_1(b_1 + b_2 + c_1 + c_2) + a_2(b_2 + 2c_1 + 3c_2) = 0$ that $\tilde{x} \in C^2([0,1]; \Re)$ is a solution of the boundary value problem (2.8), (2.9). Uniqueness follows from the fact that $y \in C^2([0,1]; \Re)$ is the unique solution of the boundary value problem $Ay = -Ag$ with $b_1 y(0) + b_2 \frac{dy}{dz}(0) = a_1 y(1) + a_2 \frac{dy}{dz}(1) = 0$.

Since $\tilde{x} \in C^2([0,1]; \Re)$, it follows that $\tilde{x} \in L_r^2([0,1])$. Since the eigenfunctions $\{\phi_n\}_{n=1}^\infty$ of the Sturm-Liouville operator $A: D \to C^0([0,1]; \Re)$ defined by (2.1), (2.2) form an orthonormal basis of $L_r^2([0,1])$, it follows that Parseval's identity holds, i.e.,

$$\|\tilde{x}\|_r^2 = \sum_{n=1}^\infty c_n^2 = \int_0^1 r(z) \tilde{x}^2(z) dz \tag{A.1}$$

where

$$c_n := \langle \phi_n, \tilde{x} \rangle = \int_0^1 r(z) \tilde{x}(z) \phi_n(z) dz, \text{ for } n = 1, 2, \ldots \tag{A.2}$$

By virtue of (2.8), (2.9) and the facts that $-\frac{1}{r(z)}\frac{d}{dz}\left(p(z)\frac{d\phi_n}{dz}(z)\right)+\frac{q(z)}{r(z)}\phi_n(z)=\lambda_n\phi_n(z)$, $b_1\phi_n(0)+b_2\frac{d\phi_n}{dz}(0)=a_1\phi_n(1)+a_2\frac{d\phi_n}{dz}(1)=0$, it follows from repeated integration by parts, that the following equalities hold for $n=1,2,\ldots$:

$$\lambda_n c_n = \int_0^1 r(z)\tilde{x}(z)\lambda_n\phi_n(z)dz = \int_0^1 r(z)\tilde{x}(z)(A\phi_n)(z)dz$$

$$= -\int_0^1 \tilde{x}(z)\frac{d}{dz}\left(p(z)\frac{d\phi_n}{dz}(z)\right)dz + \int_0^1 q(z)\tilde{x}(z)\phi_n(z)dz$$

$$= p(0)\tilde{x}(0)\frac{d\phi_n}{dz}(0) - p(1)\tilde{x}(1)\frac{d\phi_n}{dz}(1) + \int_0^1 p(z)\frac{d\tilde{x}}{dz}(z)\frac{d\phi_n}{dz}(z)dz + \int_0^1 q(z)\tilde{x}(z)\phi_n(z)dz$$

$$= p(0)\tilde{x}(0)\frac{d\phi_n}{dz}(0) - p(1)\tilde{x}(1)\frac{d\phi_n}{dz}(1) + p(1)\frac{d\tilde{x}}{dz}(1)\phi_n(1) - p(0)\frac{d\tilde{x}}{dz}(0)\phi_n(0) - \int_0^1 \phi_n(z)\frac{d}{dz}\left(p(z)\frac{d\tilde{x}}{dz}(z)\right)dz + \int_0^1 q(z)\tilde{x}(z)\phi_n(z)dz$$

$$= p(0)\left(\tilde{x}(0)\frac{d\phi_n}{dz}(0) - \frac{d\tilde{x}}{dz}(0)\phi_n(0)\right) + p(1)\left(\frac{d\tilde{x}}{dz}(1)\phi_n(1) - \tilde{x}(1)\frac{d\phi_n}{dz}(1)\right) - \int_0^1 r(z)\phi_n(z)(A\tilde{x})(z)dz =$$

$$= p(0)\left(b_1\frac{d\phi_n}{dz}(0) - b_2\phi_n(0)\right)$$

In the above equations, we have used the fact that by virtue of (2.9) and $a_1\phi_n(1)+a_2\frac{d\phi_n}{dz}(1)=0$, the homogeneous system of linear equations

$$s_1\tilde{x}(1)+s_2\frac{d\tilde{x}}{dz}(1) = 0 = s_1\phi_n(1)+s_2\frac{d\phi_n}{dz}(1)$$

has the non-zero solution

$$s_1=a_1, s_2=a_2$$

and consequently, the determinant of the matrix $\begin{bmatrix} \tilde{x}(1) & \frac{d\tilde{x}}{dz}(1) \\ \phi_n(1) & \frac{d\phi_n}{dz}(1) \end{bmatrix}$ is zero, i.e., $\tilde{x}(1)\frac{d\phi_n}{dz}(1)-\phi_n(1)\frac{d\tilde{x}}{dz}(1)=0$. Moreover, we have used the fact that $\tilde{x}(0)\frac{d\phi_n}{dz}(0)-\frac{d\tilde{x}}{dz}(0)\phi_n(0) = b_1\frac{d\phi_n}{dz}(0)-b_2\phi_n(0)$, which is a direct consequence of the facts $b_1\tilde{x}(0)+b_2\frac{d\tilde{x}}{dz}(0)=1$, $b_1\phi_n(0)+b_2\frac{d\phi_n}{dz}(0)=0$ and $b_1^2+b_2^2=1$.

Identity (2.10) is a direct consequence of (A.1) and the fact that $c_n=p(0)\left(b_1\frac{d\phi_n}{dz}(0)-b_2\phi_n(0)\right)\lambda_n^{-1}$ for $n=1,2,\ldots$. The proof is complete. ◁

**Proof of Theorem 3.1:** Since $x_0 \in X_0$, the series $\sum_{n=1}^{\infty}c_n\phi_n(z)$ with $c_n=\int_0^1 r(z)\phi_n(z)x_0(z)dz$ ($n=1,2,\ldots$) is uniformly and absolutely convergent on $[0,1]$ and satisfies $x_0(z)=\sum_{n=1}^{\infty}c_n\phi_n(z)$ for all $z\in[0,1]$. This is a direct consequence of Theorem 9.3 on page 281 in [8], Theorem 7.5.4 on page 500 in [13], the fact that $0<\lambda_1<\lambda_2<\ldots<\lambda_n<\ldots$ and the fact that every $x_0\in X_0$ satisfies $x_0(z)=\int_0^1 g(z,s)r(s)(Ax_0)(s)ds$, where $g \in C^0([0,1]^2; \Re)$ is the Green's function of the Sturm-Liouville operator $A: D \to C^0([0,1]; \Re)$ defined by (2.1), (2.2). Define:

$$\theta_n(t) := \int_0^1 r(z)\phi_n(z)f(t,z)dz, \text{ for all } t \geq 0, \ n = 1,2,... \tag{A.3}$$

Since $f \in C^1(\Re_+ \times [0,1]; \Re)$, it follows from Theorem 3.11.3.4 in [3], that the following equation holds:

$$f(t,z) = \sum_{n=1}^{\infty} \theta_n(t)\phi_n(z), \text{ for all } (t,z) \in (0,+\infty) \times (0,1) \tag{A.4}$$

Moreover, notice that the Cauchy-Schwarz inequality, in conjunction with the fact that $\|\phi_n\|_r = 1$ (for $n = 1,2,...$) and the fact that $f \in C^1(\Re_+ \times [0,1]; \Re)$, implies the following relations for all $t \geq 0$:

$$|\theta_n(t)| \leq \left(\int_0^1 r(z)|f(t,z)|^2 dz\right)^{1/2} \tag{A.5}$$

$$\dot{\theta}_n(t) = \int_0^1 r(z)\phi_n(z)\frac{\partial f}{\partial t}(t,z)dz \tag{A.6}$$

$$|\dot{\theta}_n(t)| \leq \left(\int_0^1 r(z)\left|\frac{\partial f}{\partial t}(t,z)\right|^2 dz\right)^{1/2} \tag{A.7}$$

Since (2.4) holds and since $0 < \lambda_1 < \lambda_2 < ... < \lambda_n < ...$, inequalities (A.5), (A.7) and the fact that the series $\sum_{n=1}^{\infty} c_n\phi_n(z)$ with $c_n = \int_0^1 r(z)\phi_n(z)x_0(z)dz$ ($n = 1,2,...$) is uniformly and absolutely convergent on $[0,1]$, imply that for every $T > 0$, the series

$$\sum_{n=1}^{\infty} \exp(-\lambda_n t)c_n\phi_n(z) + \sum_{n=1}^{\infty} \phi_n(z)\lambda_n^{-1}(\theta_n(t) - \theta_n(0)\exp(-\lambda_n t)) - \sum_{n=1}^{\infty} \phi_n(z)\lambda_n^{-1}\int_0^t \exp(-\lambda_n(t-s))\dot{\theta}_n(s)ds$$

is uniformly and absolutely convergent on $[0,T] \times [0,1]$. Therefore, we define $x \in C^0(\Re_+ \times [0,1]; \Re)$ by means of the formula:

$$x(t,z) := \sum_{n=1}^{\infty} \phi_n(z)\left(\exp(-\lambda_n t)c_n + \lambda_n^{-1}\theta_n(t) - \lambda_n^{-1}\theta_n(0)\exp(-\lambda_n t) - \lambda_n^{-1}\int_0^t \exp(-\lambda_n(t-s))\dot{\theta}_n(s)ds\right),$$

$$\text{for all } (t,z) \in \Re_+ \times [0,1] \tag{A.8}$$

In order, to show that the derivative $\frac{\partial x}{\partial t}(t,z)$ exists for every $(t,z) \in (0,+\infty) \times [0,1]$ and is a continuous mapping, we show that for every $0 < t_0 < T$, the series obtained (formally) by term-by-term differentiation of the right hand side of (A.8) with respect to $t$ is uniformly and absolutely convergent on $[t_0, T] \times [0,1]$. Indeed, we get from term-by-term differentiation of the right hand side of (A.8) with respect to $t$:

$$-\sum_{n=1}^{\infty} \lambda_n \exp(-\lambda_n t)c_n\phi_n(z) + \sum_{n=1}^{\infty} \phi_n(z)\theta_n(0)\exp(-\lambda_n t) + \sum_{n=1}^{\infty} \phi_n(z)\int_0^t \exp(-\lambda_n(t-s))\dot{\theta}_n(s)ds$$

Inequality (A.7) implies that $\left|\int_0^t \exp(-\lambda_n(t-s))\dot\theta_n(s)ds\right| \le \lambda_n^{-1}\max_{0\le s\le T}\left(\left|\dot\theta_n(s)\right|\right) \le \lambda_n^{-1}\left(\max_{0\le s\le T}\left(\int_0^1 r(z)\left|\frac{\partial f}{\partial t}(t,z)\right|^2 dz\right)\right)^{1/2}$, which combined with the inequalities $\lambda_n\exp(-\lambda_n t) = t^{-1}\lambda_n t\exp(-\lambda_n t) \le t_0^{-1}\exp(-1)$, $\lambda_n^2\exp(-\lambda_n t) = t^{-2}(\lambda_n t)^2\exp(-\lambda_n t) \le 4t_0^{-2}\exp(-2)$ that hold for all $t\in[t_0,T]$, the fact that $\{c_n\}_{n=1}^\infty$, $\{\theta_n(0)\}_{n=1}^\infty$ are bounded sequences (recall (A.5) and notice that since $c_n = \int_0^1 r(z)\phi_n(z)x_0(z)dz$, the Cauchy-Schwarz inequality implies $|c_n| \le \left(\int_0^1 r(z)x_0^2(z)dz\right)^{1/2}$) and (2.4), guarantees that the series is uniformly and absolutely convergent on $[t_0,T]\times[0,1]$. Therefore, $\frac{\partial x}{\partial t}(t,z)$ exists for every $(t,z)\in(0,+\infty)\times[0,1]$ and is a continuous mapping that satisfies

$$\frac{\partial x}{\partial t}(t,z) = \sum_{n=1}^\infty \phi_n(z)\left(\theta_n(0)\exp(-\lambda_n t) - \lambda_n\exp(-\lambda_n t)c_n + \int_0^t \exp(-\lambda_n(t-s))\dot\theta_n(s)ds\right) \quad (A.9)$$

for all $(t,z)\in(0,+\infty)\times[0,1]$.

Since the Green's function of the Sturm-Liouville operator $A:D\to C^0([0,1];\Re)$ defined by (2.1), (2.2), $g\in C^0([0,1]^2;\Re)$ is a $C^2$ function on each one of the triangles $0\le s\le z\le 1$ and $0\le z\le s\le 1$, having a step discontinuity in $\frac{\partial g}{\partial z}(z,s)$ on the line segment $0\le s=z\le 1$ (see Theorem 2.2 on pages 227-228 in [8]) and since $\phi_n(z) = \lambda_n\int_0^1 g(z,s)r(s)\phi_n(s)ds$ for all $z\in[0,1]$, it follows that there exists a constant $M>0$ such that

$$\max_{0\le z\le 1}\left(\left|\frac{d\phi_n}{dz}(z)\right|\right) \le M\lambda_n\max_{0\le z\le 1}\left(|\phi_n(z)|\right), \text{ for all } n=1,2,... \quad (A.10)$$

Moreover, equality (A.4) implies that the following equality holds for all $(t,z)\in(0,+\infty)\times[0,1]$:

$$\sum_{n=1}^\infty \phi_n(z)\lambda_n^{-1}\theta_n(t) = \int_0^1 r(s)g(z,s)f(t,s)ds \quad (A.11)$$

Notice that the function $z\to(y[t])(z) = \int_0^1 r(s)g(z,s)f(t,s)ds$ is simply the unique solution of the boundary value problem $(Ay[t])(z) = f(t,z)$ with $b_1(y[t])(1) + b_2\frac{d(y[t])}{dz}(0) = a_1(y[t])(1) + a_2\frac{d(y[t])}{dz}(0) = 0$ for each fixed $t\ge 0$. Consequently, the function $z\to y_t(z) = \int_0^1 r(s)g(z,s)f(t,s)ds$ is in $X_0 = \left\{x\in C^2([0,1];\Re):b_1 x(0) + b_2\frac{dx}{dz}(0) = a_1 x(1) + a_2\frac{dx}{dz}(1) = 0\right\}$ for all $t\ge 0$. Moreover, we obtain from (A.8) and (A.11):

$$x(t,z) - \int_0^1 g(z,s)r(s)f(t,s)ds = \sum_{n=1}^\infty \phi_n(z)\left(\exp(-\lambda_n t)c_n - \lambda_n^{-1}\theta_n(0)\exp(-\lambda_n t) - \lambda_n^{-1}\int_0^t \exp(-\lambda_n(t-s))\dot\theta_n(s)ds\right),$$
$$\text{for all } (t,z)\in\Re_+\times[0,1] \quad (A.12)$$

In order, to show that the derivative $\frac{\partial x}{\partial z}(t,z)$ exists for every $(t,z)\in(0,+\infty)\times[0,1]$ and is a continuous mapping, we show that for every $0<t_0<T$, the series obtained (formally) by term-by-term differentiation of the right hand side of (A.12) with respect to $z$ is uniformly and absolutely convergent on $[t_0,T]\times[0,1]$. Indeed, we get from term-by-term differentiation of the right hand side of (A.12) with respect to $z$:

$$\sum_{n=1}^{\infty}\exp(-\lambda_n t)c_n\frac{d\phi_n}{dz}(z)-\sum_{n=1}^{\infty}\frac{d\phi_n}{dz}(z)\lambda_n^{-1}\dot{p}_n(0)\exp(-\lambda_n t)-\sum_{n=1}^{\infty}\frac{d\phi_n}{dz}(z)\lambda_n^{-1}\int_0^t\exp(-\lambda_n(t-s))\ddot{p}_n(s)ds$$

Inequality (A.7) implies that $\left|\int_0^t\exp(-\lambda_n(t-s))\dot{\theta}_n(s)ds\right|\leq\lambda_n^{-1}\max_{0\leq s\leq T}\left(\left|\dot{\theta}_n(s)\right|\right)\leq\lambda_n^{-1}\max_{0\leq s\leq T}\left(\int_0^1 r(z)\left|\frac{\partial f}{\partial t}(t,z)\right|^2 dz\right)^{1/2}$, which combined with the inequalities $\lambda_n\exp(-\lambda_n t)=t^{-1}\lambda_n t\exp(-\lambda_n t)\leq t_0^{-1}\exp(-1)$, $\lambda_n^2\exp(-\lambda_n t)=t^{-2}(\lambda_n t)^2\exp(-\lambda_n t)\leq 4t_0^{-2}\exp(-2)$ that hold for all $t\in[t_0,T]$, the fact that $\{c_n\}_{n=1}^{\infty}$, $\{\theta_n(0)\}_{n=1}^{\infty}$ are bounded sequences (recall (A.5) and notice that since $c_n=\int_0^1 r(z)\phi_n(z)x_0(z)dz$, the Cauchy-Schwarz inequality implies $|c_n|\leq\left(\int_0^1 r(z)x_0^2(z)dz\right)^{1/2}$) and (2.4), (A.10), guarantees that the series is uniformly and absolutely convergent on $[t_0,T]\times[0,1]$. Therefore, $\frac{\partial x}{\partial z}(t,z)$ exists for every $(t,z)\in(0,+\infty)\times[0,1]$ and is a continuous mapping that satisfies

$$\frac{\partial x}{\partial z}(t,z)-\frac{\partial}{\partial z}\left(\int_0^1 r(s)g(z,s)f(t,s)ds\right)=\sum_{n=1}^{\infty}\exp(-\lambda_n t)c_n\frac{d\phi_n}{dz}(z)$$
$$-\sum_{n=1}^{\infty}\frac{d\phi_n}{dz}(z)\lambda_n^{-1}\theta_n(0)\exp(-\lambda_n t)-\sum_{n=1}^{\infty}\frac{d\phi_n}{dz}(z)\lambda_n^{-1}\int_0^t\exp(-\lambda_n(t-s))\dot{\theta}_n(s)ds \quad (A.13)$$

for all $(t,z)\in(0,+\infty)\times[0,1]$.

The differential equation $\frac{d^2\phi_n}{dz^2}(z)=\frac{q(z)}{p(z)}\phi_n(z)-\lambda_n\frac{r(z)}{p(z)}\phi_n(z)-\frac{1}{p(z)}\frac{dp}{dz}(z)\frac{d\phi_n}{dz}(z)$, which holds for all $z\in[0,1]$, in conjunction with the fact that $0<\lambda_1<\lambda_2<...<\lambda_n<...$ and (A.10), implies that that there exists a constant $G>0$ such that

$$\max_{0\leq z\leq 1}\left(\left|\frac{d^2\phi_n}{dz^2}(z)\right|\right)\leq G\lambda_n\max_{0\leq z\leq 1}\left(\left|\phi_n(z)\right|\right),\text{ for all }n=1,2,... \quad (A.14)$$

In order, to show that the derivative $\frac{\partial^2 x}{\partial z^2}(t,z)$ exists for every $(t,z)\in(0,+\infty)\times[0,1]$ and is a continuous mapping, we show that for every $0<t_0<T$, the series obtained (formally) by term-by-term differentiation of the right hand side of (A.13) with respect to $z$ is uniformly and absolutely convergent on $[t_0,T]\times[0,1]$. Indeed, we get from term-by-term differentiation of the right hand side of (A.13) with respect to $z$:

$$\sum_{n=1}^{\infty}\exp(-\lambda_n t)c_n \frac{d^2\phi_n}{dz^2}(z) - \sum_{n=1}^{\infty}\frac{d^2\phi_n}{dz^2}(z)\lambda_n^{-1}\theta_n(0)\exp(-\lambda_n t) - \sum_{n=1}^{\infty}\frac{d^2\phi_n}{dz^2}(z)\lambda_n^{-1}\int_0^t \exp(-\lambda_n(t-s))\dot\theta_n(s)ds$$

Inequality (A.7) implies that $\left|\int_0^t \exp(-\lambda_n(t-s))\dot\theta_n(s)ds\right| \leq \lambda_n^{-1}\max_{0\leq s\leq T}\left(|\dot\theta_n(s)|\right) \leq \lambda_n^{-1}\max_{0\leq s\leq T}\left(\int_0^1 r(z)\left|\frac{\partial f}{\partial t}(t,z)\right|^2 dz\right)^{1/2}$, which combined with the inequalities $\lambda_n \exp(-\lambda_n t) = t^{-1}\lambda_n t\exp(-\lambda_n t) \leq t_0^{-1}\exp(-1)$, $\lambda_n^2 \exp(-\lambda_n t) = t^{-2}(\lambda_n t)^2 \exp(-\lambda_n t) \leq 4t_0^{-2}\exp(-2)$ that hold for all $t\in[t_0,T]$, the fact that $\{c_n\}_{n=1}^{\infty}$, $\{\theta_n(0)\}_{n=1}^{\infty}$ are bounded sequences and (2.4), (A.14), guarantees that the series is uniformly and absolutely convergent on $[t_0,T]\times[0,1]$. Therefore, $\frac{\partial^2 x}{\partial z^2}(t,z)$ exists for every $(t,z) \in (0,+\infty)\times[0,1]$ and is a continuous mapping that satisfies

$$\frac{\partial^2 x}{\partial z^2}(t,z) - \frac{\partial^2}{\partial z^2}\left(\int_0^1 r(s)g(z,s)f(t,s)ds\right) = \sum_{n=1}^{\infty}\exp(-\lambda_n t)c_n \frac{d^2\phi_n}{dz^2}(z)$$
$$- \sum_{n=1}^{\infty}\frac{d^2\phi_n}{dz^2}(z)\lambda_n^{-1}\theta_n(0)\exp(-\lambda_n t) - \sum_{n=1}^{\infty}\frac{d^2\phi_n}{dz^2}(z)\lambda_n^{-1}\int_0^t \exp(-\lambda_n(t-s))\dot\theta_n(s)ds \quad (A.15)$$

for all $(t,z) \in (0,+\infty)\times[0,1]$.

It follows that the mapping $x$ defined by (A.8) is of class $C^0(\Re_+ \times [0,1];\Re) \cap C^1((0,+\infty)\times[0,1];\Re)$ satisfying $x(t,\cdot) \in X_0$ for all $t\geq 0$ and $x(0,z) = x_0(z)$ for all $z\in[0,1]$. Equation (3.1) is a direct consequence of (A.9), (A.12), the fact that $A\phi_n = \lambda_n\phi_n$ and the fact that the function $z \to (y[t])(z) = \int_0^1 r(s)g(z,s)f(t,s)ds$ is simply the unique solution of the boundary value problem $(Ay[t])(z) = f(t,z)$ with $b_1(y[t])(1) + b_2\frac{d(y[t])}{dz}(0) = a_1(y[t])(1) + a_2\frac{d(y[t])}{dz}(0) = 0$ for each fixed $t\geq 0$.

Finally, uniqueness follows from a standard argument, which is described next. Consider two functions $x_i \in C^0(\Re_+ \times[0,1];\Re)\cap C^1((0,+\infty)\times[0,1];\Re)$ ($i=1,2$) satisfying $x_i[t] \in X_0$ for all $t\geq 0$, $x_i(0,z) = x_0(z)$ for all $z\in[0,1]$ and

$$\frac{\partial x_i}{\partial t}(t,z) - \frac{1}{r(z)}\frac{\partial}{\partial z}\left(p(z)\frac{\partial x_i}{\partial z}(t,z)\right) + \frac{q(z)}{r(z)}x_i(t,z) = f(t,z), \text{ for all } (t,z)\in(0,+\infty)\times(0,1),\ i=1,2 \quad (A.16)$$

It follows that the function $\bar x = x_1 - x_2$ is of class $C^0(\Re_+\times[0,1];\Re)\cap C^1((0,+\infty)\times[0,1];\Re)$ satisfying $\bar x[t] \in X_0$ for all $t\geq 0$, $\bar x(0,z) = 0$ for all $z\in[0,1]$ and

$$\frac{\partial \bar x}{\partial t}(t,z) - \frac{1}{r(z)}\frac{\partial}{\partial z}\left(p(z)\frac{\partial \bar x}{\partial z}(t,z)\right) + \frac{q(z)}{r(z)}\bar x(t,z) = 0, \text{ for all } (t,z)\in(0,+\infty)\times(0,1) \quad (A.17)$$

Since the eigenfunctions $\{\phi_n\}_{n=1}^{\infty}$ of the Sturm-Liouville operator $A:D\to C^0([0,1];\Re)$ defined by (2.1), (2.2) form an orthonormal basis of $L_r^2([0,1])$, it follows that Parseval's identity holds, i.e.,

$$\|\bar x[t]\|_r^2 = \sum_{n=1}^{\infty}\xi_n^2(t), \text{ for all } t\geq 0 \quad (A.18)$$

where

$$\xi_n(t) := \langle \phi_n, \bar{x}[t] \rangle = \int_0^1 r(z)\bar{x}(t,z)\phi_n(z)dz, \text{ for } n = 1,2,... \qquad (A.19)$$

By virtue of (2.5) and the facts that $\bar{x}[t] \in X_0$, $-\frac{1}{r(z)}\frac{d}{dz}\left(p(z)\frac{d\phi_n}{dz}(z)\right) + \frac{q(z)}{r(z)}\phi_n(z) = \lambda_n \phi_n(z)$, $b_1\phi_n(0) + b_2\frac{d\phi_n}{dz}(0) = a_1\phi_n(1) + a_2\frac{d\phi_n}{dz}(1) = 0$, it follows from repeated integration by parts, that the following equalities hold for all $t > 0$:

$$\dot{\xi}_n(t) = \int_0^1 r(z)\frac{\partial \bar{x}}{\partial t}(t,z)\phi_n(z)dz$$

$$= \int_0^1 \frac{\partial}{\partial z}\left(p(z)\frac{\partial \bar{x}}{\partial z}(t,z)\right)\phi_n(z)dz - \int_0^1 q(z)\bar{x}(t,z)\phi_n(z)dz$$

$$= p(1)\frac{\partial \bar{x}}{\partial z}(t,1)\phi_n(1) - p(0)\frac{\partial \bar{x}}{\partial z}(t,0)\phi_n(0) - \int_0^1 p(z)\frac{\partial \bar{x}}{\partial z}(t,z)\frac{d\phi_n}{dz}(z)dz - \int_0^1 q(z)\bar{x}(t,z)\phi_n(z)dz =$$

$$= p(1)\left(\frac{\partial \bar{x}}{\partial z}(t,1)\phi_n(1) - \bar{x}(t,1)\frac{d\phi_n}{dz}(1)\right) + p(0)\left(\frac{d\phi_n}{dz}(0)\bar{x}(t,0) - \frac{\partial \bar{x}}{\partial z}(t,0)\phi_n(0)\right) + \int_0^1 \bar{x}(t,z)\left[\frac{d}{dz}\left(p(z)\frac{d\phi_n}{dz}(z)\right) - q(z)\phi_n(z)\right]dz$$

$$= p(1)\left(\frac{\partial \bar{x}}{\partial z}(t,1)\phi_n(1) - \bar{x}(t,1)\frac{d\phi_n}{dz}(1)\right) + p(0)\left(\frac{d\phi_n}{dz}(0)\bar{x}(t,0) - \frac{\partial \bar{x}}{\partial z}(t,0)\phi_n(0)\right) - \int_0^1 r(z)\bar{x}(t,z)(A\phi_n)(z)dz = -\lambda_n \xi_n(t)$$

In the above equations, we have used the fact that by virtue of (2.5) and the facts $\bar{x}[t] \in X_0$, $b_1\phi_n(0) + b_2\frac{d\phi_n}{dz}(0) = a_1\phi_n(1) + a_2\frac{d\phi_n}{dz}(1) = 0$, the homogeneous systems of linear equations

$$\begin{array}{cc} s_1 \bar{x}(t,1) + s_2 \frac{\partial \bar{x}}{\partial z}(t,1) = 0 & s_1 \bar{x}(t,0) + s_2 \frac{\partial \bar{x}}{\partial z}(t,0) = 0 \\ & \text{and} \\ s_1 \phi_n(1) + s_2 \frac{d\phi_n}{dz}(1) = 0 & s_1 \phi_n(0) + s_2 \frac{d\phi_n}{dz}(0) = 0 \end{array}$$

have the non-zero solutions

$$s_1 = a_1, \quad s_2 = a_2 \text{ and } s_1 = b_1, \quad s_2 = b_2$$

and consequently, the determinants of the matrices $\begin{bmatrix} \bar{x}(t,1) & \frac{\partial \bar{x}}{\partial z}(t,1) \\ \phi_n(1) & \frac{d\phi_n}{dz}(1) \end{bmatrix}$, $\begin{bmatrix} \bar{x}(t,0) & \frac{\partial \bar{x}}{\partial z}(t,0) \\ \phi_n(0) & \frac{d\phi_n}{dz}(0) \end{bmatrix}$ are zero, i.e.,

$$\bar{x}(t,1)\frac{d\phi_n}{dz}(1) - \phi_n(1)\frac{\partial \bar{x}}{\partial z}(t,1) = 0 \text{ and } \bar{x}(t,0)\frac{d\phi_n}{dz}(0) - \phi_n(0)\frac{\partial \bar{x}}{\partial z}(t,0) = 0.$$

Integrating the above differential equations, we obtain for all $0 < T \leq t$ and $n = 1,2,...$:

$$\xi_n(t) = \exp(-\lambda_n(t-T))\xi_n(T) \qquad (A.20)$$

Continuity of the mapping $\Re_+ \ni T \to \xi_n(T)$ and (A.19), (A.20) in conjunction with the fact that $\bar{x}(0,z) = 0$ for all $z \in [0,1]$ implies that $\xi_n(t) = 0$ for all $t \geq 0$ and $n = 1,2,...$. We conclude from (A.18) that $\bar{x}[t] \equiv 0$, for all $t \geq 0$. Consequently, the two functions $x_i \in C^0(\Re_+ \times [0,1]; \Re) \cap C^1((0,+\infty) \times [0,1]; \Re)$ ($i = 1,2$) satisfying $x_i[t] \in X_0$ for all $t \geq 0$, $x_i(0,z) = x_0(z)$ for all $z \in [0,1]$ are identical.

The proof is complete. ◁